\documentclass [12pt] {amsart}
\usepackage{enumerate}
\usepackage{palatino}
\usepackage{latexsym,amssymb, amsmath}
\usepackage{graphicx}

\usepackage{simplemargins}
\setleftmargin{2.75cm}
\setrightmargin{2.75cm}
\settopmargin{2.0in}
\setbottommargin{2.0in}
\setallmargins{2.75cm}

\def\pf{\noindent\emph{Proof: }}
\def\stop{\hfill$\Box$}

\usepackage{amsmath, amsfonts, amssymb, amsthm}
\newtheorem{thm}{Theorem}
\newtheorem{cor}[thm]{Corollary}
\newtheorem{lemma}[thm]{Lemma}

\newtheorem{defn}[thm]{Definition}
\newtheorem{prop}[thm]{Proposition}

\numberwithin{thm}{section}

\DeclareMathOperator{\Hess}{Hess}
\DeclareMathOperator{\diver}{div}

\begin{document}
	
	\title[Twice Differentiable Functions with Continuous Laplacian and Bounded Hessian]{Construction of Twice Differentiable Functions with Continuous Laplacian and Bounded Hessian}

	\author{Yifei Pan}
	\address{Department of Mathematical Sciences\\	Purdue University Fort Wayne\\	Fort Wayne, Indiana 46805}
	\email{pan@pfw.edu}

\author{ Yu Yan }
\address {Department of Mathematics and Computer Science\\ Biola University\\La Mirada, California 90639}
\email {yu.yan@biola.edu}
	
	\begin{abstract}
	We construct examples of twice differentiable functions in $\mathbb{R}^n$ with continuous Laplacian and bounded Hessian.	The same construction is also applicable to higher order differentiability, the Monge-Amp\`ere equation, and the mean curvature equation for hypersurfaces. 
	\end{abstract}
	
	
	\maketitle
	\newtheorem{Thm}{Theorem}
	\newtheorem{Lemm}{Lemma}
	\newtheorem{Cor}{Corollary}

\section{Introduction}
\label{section:introduction}

\vspace{.1in}
  One of the central themes in elliptic theory is to understand the extent to which the Laplacian controls a function's regularity.  By the classical Schauder theory, if $\Delta u \in C^{0, \alpha} \, (0 < \alpha <1)$, then $ u \in C^{2, \alpha}$, but this is false when $\alpha =0$. There are well known functions (\cite{Oton}) that have continuous Laplacian yet are not $C^2$; for example,  
  \begin{equation*}
  	w(x,y)= \left\{
  	\begin{array}{l l }      
  		(x^2-y^2)\ln (-\ln (x^2+y^2)) & \hspace{.4in}  0<x^2+y^2\leq \frac{1}{4} , \\
  		\noalign{\smallskip}
  		0 & \hspace{.4in}  (x,y)=(0,0).	
  	\end{array}\right.
  \end{equation*}
  
  \noindent
  Precisely, this function is not $C^2$ because it fails to be twice differentiable at the origin.  This phenomenon raises the natural question of whether continuous Laplacian and twice differentiability everywhere would be sufficient to guarantee continuous twice differentiability, i.e. $C^2$. To answer this question, we constructed (\cite{Pan-Yan_unbounded}) a large family of functions that are twice differentiable everywhere with continuous Laplacian but unbounded Hessian, hence not $C^2$. In light of the existence of such functions one naturally wonders, \textit{is an everywhere twice differentiable function with continuous Laplacian and bounded Hessian necessarily $C^2$}?
  
  \vspace{.05in}
  If twice differentiability everywhere is not required, then there are already simple examples of functions with continuous Laplacian and bounded Hessian that are twice differentiable except at the origin, such as (\cite{Oton})
  \begin{equation*}
  	\phi(x, y) =	\left\{
  	\begin{array}{l l }      
  		(x^2-y^2)\sin (\ln (-\ln (x^2+y^2))) & \hspace{.4in}  0<x^2+y^2\leq \frac{1}{4} , \\
  		\noalign{\smallskip}
  		0 & \hspace{.4in}  (x,y)=(0,0).
  	\end{array}\right.
  \end{equation*}

  \vspace{.05in}
  Our first observation is that any possible counter-examples must be non-radially symmetric because of the following result.
  \begin{prop}
  	\label{prop:radially_symmetric}
  	Any twice differentiable, radially symmetric function with continuous Laplacian must be $C^2$.
  \end{prop}

\vspace{.1in}

  A similar situation (of twice differentiability) arose in the study of parabolic equations related to the mean curvature flow.  
   If for each $s$  the level sets $$ M_t =  \{ x \,\, | \,\, v(x,t) =s, \,\, \text{where} \,\, v: \mathbb{R}^{n+1} \times \mathbb{R} \to \mathbb{R}\}$$ evolve by the mean curvature flow, then $v$ satisfies the level set equation
  \begin{equation*}
  	\partial _t v = |\triangledown v| \diver \left ( \frac{\triangledown v}{|\triangledown v|} \right ).
  \end{equation*}  
In general, $v$ is a weak solution in the viscosity sense, so it may not be differentiable.  When the initial hypersurface is mean convex, Evans and Spruck (\cite{Evans_Spruck}) showed that  $v(x,t)=u(x)-t$, where $u$ is a Lipschitz function and satisfies (in the viscosity sense)
 \begin{equation*}
	-1 = |\triangledown u| \diver \left ( \frac{\triangledown u}{|\triangledown u|} \right ).
\end{equation*} 
The function $u$ is called the arrival time. Colding and Minicozzi proved in \cite{Colding_Min_2} that surprisingly $u$ is twice differentiable everywhere with bounded Hessian, furthermore, it satisfies the equation everywhere in the classical sense. It is intriguing to know whether $\Delta u$ is continuous without $u$ being $C^2$.  If true, that would provide ample examples of twice differentiable functions with continuous Laplacian and bounded Hessian that are not $C^2$.  
 


\vspace{.05in}

\vspace{.05in}
In this paper we will show the existence of twice differentiable functions with continuous Laplacian and bounded (but discontinuous) Hessian, thus answering the aforementioned question in the negative.  That is, an everywhere twice differentiable function with continuous Laplacian and bounded Hessian is not necessarily $C^2$.

\vspace{.05in}
\begin{thm}
	\label{thm:bounded_2nd_order}
	Given any $C^2$ function $\varphi: (0, \infty) \to \mathbb{R}$ satisfying
	\begin{equation*}
		\lim _{s \to \infty} \varphi  (s) = \infty,  \hspace{.2in} \lim _{s \to \infty} \varphi ' (s) = 0,  \hspace{.2in} \text{and} \hspace{.2in} \lim _{s \to \infty} \varphi ''(s) = 0,
		\end{equation*}
	
	\noindent
	there is a compactly supported function $u: \mathbb{R}^n \to \mathbb{R}$ depending on $\varphi$, such that it is twice differentiable everywhere in $\mathbb{R}^n$ with continuous Laplacian and bounded Hessian, but $u$ is not in $C^2(\mathbb{R}^n)$.
\end{thm}

\vspace{.05in}
\noindent
Some examples of a function $\varphi$ satisfying the given conditions are  $\varphi(s)= s^{\alpha}$ with $ 0<\alpha <1$, $\varphi(s)=\ln (s)$, or $\varphi(s)= \ln   \ln \cdots  \ln s  $ where $s>e^{e^{\cdots ^ e}}$.

\vspace{.05in}

\vspace{.05in}
By imposing conditions on higher order derivatives of $\varphi$, one can generalize Theorem \ref{thm:bounded_2nd_order} to the higher order case.

\vspace{.05in}
\begin{thm}
	\label{thm: bounded_kth_order}
	Given any $C^{k+2}$ function $\varphi: (0, \infty) \to \mathbb{R}$ satisfying
	\begin{equation*}
		\lim _{s \to \infty} \varphi  (s) = \infty  \hspace{.2in} \text{and} \hspace{.2in} \lim _{s \to \infty} \varphi ' (s) = \cdots = \lim _{s \to \infty} \varphi ^{(k+2)}(s) = 0,
	\end{equation*}
	there is a compactly supported function $u:\mathbb{R}^n \to \mathbb{R}$ depending on $\varphi$, such that $u$ is $(k+2)$-times differentiable everywhere in $\mathbb{R}^n$, $\Delta u$ is $C^k$, $D^{k+2}u$ is bounded, but $u$ is not in $C^{k+2}(\mathbb{R}^n)$.
\end{thm}

\vspace{.05in}

\vspace{.1in}
We can also construct examples for the fully nonlinear Monge-Amp\`ere equations by assuming additional conditions on $\varphi$.

\begin{thm}
	\label{thm: monge-ampere}
			Given any $C^2$ function $\varphi: (0, \infty) \to \mathbb{R}$ satisfying
		\begin{equation*}
			\lim _{s \to \infty} \varphi  (s) = \infty,  \hspace{.2in} \lim _{s \to \infty} \varphi ' (s) = 0,  \hspace{.2in} \lim _{s \to \infty} \varphi ''(s) = 0,
		\end{equation*}
	and
	\begin{equation*}
		\lim_{s \to \infty} \varphi ^{n-1} (s) \varphi '(s) =0,  \hspace{.4in} 	\lim_{s \to \infty} \varphi ^{n-1} (s) \varphi ''(s) =0  ,
	\end{equation*}
		
		\vspace{.05in}
		\noindent
		there is a compactly supported function $u: \mathbb{R}^n \to \mathbb{R}$ $(n \geq 3)$ depending on $\varphi$, such that it is twice differentiable everywhere in $\mathbb{R}^n$ and $\det D^2 u$ is continuous, but $u$ is not in $C^2(\mathbb{R}^n)$. Here we can construct $u$ to have either bounded or unbounded Hessian.
	\end{thm}

\vspace{.05in}
\noindent
Some choices of such $\varphi$ are $\varphi(s)=\ln s$ and $\varphi(s)=\ln \cdots \ln s$ $(s>e^{e^{\cdots ^ e}})$.  We would like to also point out that the dimension $n=2$ case is still open.

\vspace{.1in}
One geometric application of Theorem \ref{thm:bounded_2nd_order} is the existence of certain hypersurfaces in $\mathbb{R}^{n+1}.$  

\begin{cor}
	\label{cor:hypersurface}
	There exists a twice differentiable hypersurface in $\mathbb{R}^{n+1}$ that is not $C^2$, but its mean curvature is continuous, and its second fundamental form is bounded.
\end{cor}

\vspace{.1in}
For related works in this direction we refer the reader to \cite{Coffman-Pan-Zhang}, \cite{Pan-Yan_unbounded}, and \cite{Coffman-Pan}.  The rest of the paper is organized as follows. In section \ref{section:2nd_proof} we prove Theorem \ref{thm:bounded_2nd_order} and outline the proof of Theorem \ref{thm: bounded_kth_order}.  The proof of Theorem \ref{thm: monge-ampere} and Corollary \ref{cor:hypersurface} are given in Section \ref{section:monge-ampere} and Section \ref{section:geometric}, respectively.  Finally we prove Proposition \ref{prop:radially_symmetric} in Section \ref{section:radial}.

\vspace{.05in}

	\vspace{.2in}

\section{Examples for the Laplacian Operator}
\label{section:2nd_proof}

\vspace{.1in}
\noindent
Recall that  $\varphi$ is a function satisfying
	\begin{equation}
	\label{eqn: condition_phi}
	\lim _{s \to \infty} \varphi  (s) = \infty,  \hspace{.2in} \lim _{s \to \infty} \varphi ' (s) = 0,  \hspace{.2in} \text{and} \hspace{.2in} \lim _{s \to \infty} \varphi ''(s) = 0.
\end{equation}

\noindent
Thus for $x=(x_1,...,x_n) \in \mathbb{R}^n$,  $\displaystyle \lim_{|x| \to 0} \varphi(-\ln |x|^2) = \infty$. Nevertheless, the product of $\varphi(-\ln |x|^2)$ with positive powers of $|x|$ is well controlled, and this fact plays a key role in our construction. The following two simple lemmas, which were proved in \cite{Pan-Yan_unbounded}, will be used repeatedly.  

\begin{lemma}
	\label{lemma:x-times-varphi}
	(Lemma 2.1 in \cite{Pan-Yan_unbounded}) For any  $\beta >0$ and $\varphi$ satisfying (\ref{eqn: condition_phi}),
	\begin{equation*}
		\label{eqn: z_times_varphi}
		\lim _{|x|  \to 0}  |x|^\beta  \varphi \left (-\ln  |x|^2 \right) =0.
	\end{equation*} 
	
\end{lemma}

\pf See page 4 of \cite{Pan-Yan_unbounded}.

\stop

\begin{lemma}
	\label{lemma:t_times_x^t_time_varphi}
(Lemma 2.2 in \cite{Pan-Yan_unbounded})	For any $0 < \beta \leq 1$ and $\varphi$ satisfying 	(\ref{eqn: condition_phi}), there is a constant $ C_{\varphi}$ depending only on $\varphi$, such that  $$   \sup_{\substack{|x| \leq \frac{2}{3} \\ 0<\beta \leq 1}} \beta |x|^\beta \left \vert  \varphi \left (-\ln  |x|^2 \right) \right \vert \leq  C_{ \varphi}. $$
\end{lemma}

\pf See pages 4-5 of \cite{Pan-Yan_unbounded}.

\stop

\vspace{.1in}
\noindent
The next lemma is a simple observation that if we remove the coefficient $\beta$ from the function in Lemma \ref{lemma:t_times_x^t_time_varphi}, then the remaining function  $ |x|^{\beta} \varphi \left (-\ln  |x|^2 \right)     $ will not be uniformly bounded as $\beta \to 0$. This will be crucial to our proofs in the later sections. 

\begin{lemma}
	\label{lem: x-time-varphi-infinity}
	For any $\varphi$ satisfying (\ref{eqn: condition_phi}), we have
	$$ \sup_{|x| \leq \frac{1}{2}}  |x|^\beta   \varphi \left (-\ln  |x|^2 \right)  \to \infty \hspace{.2in} \text{as} \hspace{.2in} \beta \to 0. $$
\end{lemma}

\pf For any $0 < \beta <1$, choose a point $x_{\beta} \in \mathbb{R}^n$ such that $\displaystyle |x_{\beta}|=e^{-\frac{1}{\beta}}$,
 then $$|x_{\beta}|^{\beta}   \varphi \left ( -\ln  |x_{\beta}|^2  \right ) = e^{-1}   \varphi \left ( \frac{2}{\beta}  \right ) \to \infty \hspace{.2in} \text{as} \hspace{.1in} \beta \to 0 . $$ 
When $\beta$ is sufficiently small, $|x_{\beta}|= e^{-\frac{1}{\beta}} < \frac{1}{2}$, so $$|x_{\beta}|^{\beta}   \varphi \left ( -\ln  |x_{\beta}|^2  \right ) \leq \sup_{|x| \leq \frac{1}{2}}  |x|^\beta   \varphi \left (-\ln  |x|^2 \right) .$$
Therefore $$ \sup_{|x| \leq \frac{1}{2}}  |x|^\beta   \varphi \left (-\ln  |x|^2 \right)  \to \infty \hspace{.2in} \text{as} \hspace{.2in} \beta \to 0. $$

\stop

\vspace{.1in}
The first two steps in the proof of Theorem \ref{thm:bounded_2nd_order} is the same as that in the proof of Theorem 1.1 in \cite{Pan-Yan_unbounded}. It was proved in \cite{Pan-Yan_unbounded} that for any $|x| \leq \frac{1}{2}$, the function $v(x)$ defined by
\begin{equation*}
	\label{eqn:defn_v}
	v(x) = 	\left\{
	\begin{array}{l l }      
			x_1x_2 \varphi \left (-\ln |x|^2 \right) & \hspace{.1in} 0 < |x| \leq  \frac{1}{2}, \\
			\noalign{\medskip}
				0 &   \hspace{.1in} x=0
	\end{array}\right.
\end{equation*}
has continuous Laplacian and unbounded Hessian, but it is not twice differentiable at $0$. This is
the building block function of the subsequent construction. We will then modify $v$ into a one-parameter family of functions $u_t$
defined below that is $C^2$ everywhere, including at the origin. In the last step, we will combine
an appropriate sequence of such $u_t$ into a function $u$ that loses the $C^2$ regularity but retains twice differentiability at the origin, furthermore, this function $u$ has continuous Laplacian and bounded Hessian.


\begin{defn}
	\label{defn: u_t}
	Let $\eta: [0, \infty) \to [0,1]$ be a fixed, non-increasing $C^{\infty}$ function such that
	\begin{equation}
		\label{eqn: condition_eta}
		\eta(s) \equiv 1 \,\, \text {for } \,\, 0 \leq s \leq \frac{1}{2}  \hspace{.2in}   \text{and} \hspace{.2in} \eta(s) \equiv 0 \,\, \text {for } \,\, s \geq \frac{2}{3}  . 
	\end{equation}
	
\vspace{.05in}
	\noindent 
	For any $0<t\leq \frac{1}{2}$, define a function $u_t: \mathbb{R}^n \to \mathbb{R} \,\, (n \geq 2)$  by
	\begin{equation}
		\label{eqn: defn_u_t}
		u_t(x) = 	\left\{
		\begin{array}{l l l}      
			0 &   \hspace{.1in} x=0, \\
			\noalign{\medskip}
			\eta(|x|)x_1 x_2|x|^{2t}\varphi  (-\ln  |x|^2  ) & \hspace{.1in} 0 < |x|< 1, \\
			\noalign{\medskip}
			0 &	\hspace{.1in} |x| \geq 1.
		\end{array}\right.
	\end{equation}
	
\end{defn}

\vspace{.1in}

The following lemma was proved in \cite{Pan-Yan_unbounded} and will be essential to our construction. 

\begin{lemma}
	\label{lemma:u_t_derivative_bounded} 
(Lemma 3.2 of \cite{Pan-Yan_unbounded})	There is a constant $ C_{\eta, \varphi}$ depending only on $\eta$ and $\varphi$, such that 

\begin{equation}
	\label{eqn:u_t_bounded}
	\sup _{x \in \mathbb{R}^n} \left \vert  u_t (x) \right \vert \leq C_{\eta, \varphi} . \hspace{1.8in}
\end{equation}

\begin{equation}
	\label{eqn:u_t_derivative_bounded}
	\sup _{x \in \mathbb{R}^n} \left \vert \frac{\partial  u_t}{\partial x_j} (x) \right \vert \leq C_{\eta, \varphi} \hspace{.4in} \text{for} \hspace{.2in} j=1, ..., n.
\end{equation}

\begin{equation}
	\label{eqn:u_t_derivative_bounded_second_diagonal}
	\sup _{x \in \mathbb{R}^n} \left \vert \frac{\partial ^2 u_t}{\partial x_j^2} (x) \right \vert \leq C_{\eta, \varphi} \hspace{.4in} \text{for} \hspace{.2in} j=1, ..., n.
\end{equation}

\end{lemma}

\vspace{.1in}

\pf See pages 8-11 of \cite{Pan-Yan_unbounded}.
\stop

\vspace{.1in}

In addition, we will also need to use the fact that except $\frac{\partial ^2 u_t}{\partial x_1x_2}$, all the other non-diagonal entries of $D^2 u_t$ are also uniformly bounded.

 \begin{lemma}
 	\label{lemma:u_t_derivative_bounded_non-diagonal} 
 	There is a constant $ C_{\eta, \varphi}$ depending only on $\eta$ and $\varphi$, such that 	
 	\begin{equation}
 		\label{eqn:u_t_derivative_bounded_x_1_x_j}
 		\sup _{x \in \mathbb{R}^n} \left \vert \frac{\partial ^2 u_t}{\partial x_1 \partial x_j} (x) \right \vert \leq  C_{\eta, \varphi} \hspace{.4in} \text{for} \hspace{.2in} j=3, ..., n.
 	\end{equation}
 	
 		\begin{equation}
 		\label{eqn:u_t_derivative_bounded_x_2_x_j}
 		\sup _{x \in \mathbb{R}^n} \left \vert \frac{\partial ^2 u_t}{\partial x_2 \partial x_j} (x) \right \vert \leq  C_{\eta, \varphi} \hspace{.4in} \text{for} \hspace{.2in} j=3, ..., n.
 	\end{equation}
 \noindent
 and
 	\begin{equation}
 		\label{eqn:u_t_derivative_bounded_x_i_x_j}
 		\sup _{x \in \mathbb{R}^n} \left \vert \frac{\partial ^2 u_t}{\partial x_i \partial x_j} (x) \right \vert \leq  C_{\eta, \varphi} \hspace{.4in} \text{for} \hspace{.2in} i,j=3, ..., n \,\, \text{and} \,\, i \neq j.
 	\end{equation}
 	
 \end{lemma}

\vspace{.1in}
\pf By the definition of $u_t$, $$ \sup _{x \in \mathbb{R}^n} \left \vert \frac{\partial ^2 u_t}{\partial x_1 \partial x_j} (x) \right \vert =  \sup _{|x| \leq \frac{2}{3} } \left \vert  \frac{\partial ^2 u_t}{\partial x_1 \partial x_j} (x) \right \vert , \hspace{.4in}  \sup _{x \in \mathbb{R}^n} \left \vert \frac{\partial ^2 u_t}{\partial x_2 \partial x_j } (x) \right \vert =  \sup _{|x| \leq \frac{2}{3}} \left \vert \frac{\partial ^2 u_t}{\partial x_2 \partial x_j} (x) \right \vert , \hspace{.2in} $$ and   $$ \sup _{x \in \mathbb{R}^n} \left \vert \frac{\partial ^2 u_t}{ \partial x_i \partial x_j } (x) \right \vert  =  \sup _{|x| \leq \frac{2}{3}} \left \vert \frac{\partial ^2 u_t}{\partial x_i \partial x_j} (x) \right \vert, $$ so we assume $|x| \leq \frac{2}{3}$.  Furthermore, since $0<t \leq \frac{1}{2}$, when $\frac{1}{2} \leq |x| \leq \frac{2}{3}$, we know $\left \vert \frac{\partial ^2 u_t}{\partial x_1 \partial x_j} (x)  \right \vert $, $\left \vert \frac{\partial ^2 u_t}{\partial x_2 \partial x_j}  (x) \right \vert $, and $\left \vert \frac{\partial ^2 u_t}{ \partial x_i \partial x_j } (x) \right \vert$ are all bounded by a constant depending on $\eta$ and $\varphi$ and independent of $t$.  Therefore it remains to show that when $|x| < \frac{1}{2}$, $\left \vert \frac{\partial ^2 u_t}{\partial x_1 \partial x_j} (x) \right \vert $, $\left \vert \frac{\partial ^2 u_t}{\partial x_2 \partial x_j } (x)  \right \vert $, and $\left \vert \frac{\partial ^2 u_t}{ \partial x_i \partial x_j } (x) \right \vert$ are all bounded by a constant independent of $t$.

\vspace{.1in}
\noindent
We first prove (\ref{eqn:u_t_derivative_bounded_x_1_x_j}).  When  $|x| \leq \frac{1}{2}$,  by definition $\eta(|x|) \equiv 1$, so for $j= 3,...,n$,
\allowdisplaybreaks
\begin{eqnarray}
	\label{eqn:2nd_derivatives_u_t_x_1_x_j}
	\hspace{.5in} 	\frac{\partial ^2 u_t}{\partial x_1 \partial x_j} (x) 
	& = &   2t  x_2 x_j   |x|^{2t-2} \varphi (-\ln  |x|^2  )  +  2t(2t-2)   x_1 ^2 x_2x_j |x|^{2t-4} \varphi (-\ln  |x|^2  )    \\	
	& + &  (4-8t)  x_1^2 x_2 x_j |x|^{2t-4} \varphi '  (-\ln |x|^2  )  + 4  x_1^2 x_2 x_j  |x|^{2t-4} \varphi '' (-\ln |x|^2  ).  \nonumber
\end{eqnarray}

\vspace{.05in}
\noindent
  The third term, $$ (4-8t) x_1^2 x_2 x_j |x|^{2t-4} \varphi '  (-\ln  |x|^2  ) ,$$
is bounded by  $$C|x|^{2t} \left \vert \varphi '  (-\ln |x|^2  ) \right \vert ,$$  which is further bounded by $$C \left \vert \varphi '  (-\ln |x|^2  ) \right \vert$$ since $|x| < 1$.  Because of (\ref{eqn: condition_phi}), we know $ \varphi '  (-\ln |x|^2  ) $ has a removable discontinuity at $0$, therefore on the closed set $|x| \leq \frac{1}{2}$ it is bounded by a constant depending only on $\varphi$.  Thus the third term is bounded by a constant depending only on $\varphi$.  Similarly, the fourth term is also bounded by a constant depending only on $\varphi$.

\vspace{.1in}
\noindent
The first term, $$2t x_2 x_j |x|^{2t-2} \varphi  (-\ln  |x|^2  ),$$ and the second term, $$ 2t(2t-2) x_1^2 x_2 x_j |x|^{2t-4} \varphi  (-\ln |x|^2  ),$$ are bounded by $$Ct|x|^{2t}\left \vert  \varphi  (-\ln |x|^2  ) \right \vert.  $$ By Lemma \ref{lemma:t_times_x^t_time_varphi}, $$  2t|x|^{2t}\left \vert  \varphi  (-\ln |x|^2  ) \right \vert \leq   C_{\varphi}, $$

\vspace{.05in}
\noindent
 where $ C_{\varphi}$ depends only on $\varphi$. Hence the first and second terms are bounded by a constant depending on $\varphi$.  Therefore, we have proved that all the terms in (\ref{eqn:2nd_derivatives_u_t_x_1_x_j}) are uniformly bounded by a constant independent of $t$.  The bound for $ \frac{\partial ^2 u_t}{\partial x_2 x_j} $  in (\ref{eqn:u_t_derivative_bounded_x_2_x_j}) can be established in the same way.

\vspace{.05in}
\noindent
It remains to show that $\frac{\partial ^2 u_t}{\partial x_i \partial x_j}$ is uniformly bounded, where $i,j=3,...,n$ and $i \neq j$.  Again since $\eta(|x|) \equiv 1$ for $|x| \leq \frac{1}{2}$,
\allowdisplaybreaks
\begin{eqnarray}
	\label{eqn:2nd_derivatives_u_t_x_i_x_j}
	\hspace{.5in} 	\frac{\partial ^2 u_t}{\partial x_i \partial x_j} (x) 
	& = &      2t(2t-2)  x_1  x_2x_i x_j |x|^{2t-4} \varphi (-\ln  |x|^2  )   \\	
	& + &  (4-8t)   x_1 x_2 x_i x_j |x|^{2t-4} \varphi '  (-\ln |x|^2  )  + 4    x_1 x_2 x_i x_j  |x|^{2t-4} \varphi '' (-\ln |x|^2  ). \nonumber
\end{eqnarray}

\vspace{.05in}
\noindent
  The first term in (\ref{eqn:2nd_derivatives_u_t_x_i_x_j}) is bounded by $Ct|x|^{2t} \left \vert \varphi (-\ln |x|^2) \right \vert $. In the analysis of (\ref{eqn:2nd_derivatives_u_t_x_1_x_j}) we have shown that such a function is bounded by a constant independent of $t$. The second and third terms in (\ref{eqn:2nd_derivatives_u_t_x_i_x_j}) are bounded by   $C \left \vert \varphi ' (-\ln |x|^2) \right \vert$ and $ C \left \vert \varphi '' (-\ln |x|^2) \right \vert $, respectively. Since $\displaystyle \lim _{s \to 0} \varphi ' (s) = \lim _{s \to 0} \varphi '' (s) =0, $  these two functions are also bounded by a constant independent of $t$. Therefore, $\frac{\partial ^2 u_t}{\partial x_i \partial x_j} (x)$ is uniformly bounded, and this proves (\ref{eqn:u_t_derivative_bounded_x_i_x_j}).

\stop

\vspace{.15in}

 Now we are ready to construct the function $u$ in Theorem \ref{thm:bounded_2nd_order} by ``piecing together" a sequence of functions $u_{t_k}$ as follows.

\vspace{.05in}
\noindent
Choose two decreasing sequences of numbers $R_k \to 0$ and $r_k \to 0$, such that $$ R_k > r_k, $$ and for geometric reasons that will be explained later we also require
\begin{equation}
	\label{eqn:R_k_r_k}
	R_k -r_k > R_{k+1}+r_{k+1};
\end{equation}
for example, we may choose $R_k=10^{-k}$ and $r_k=10^{-(k+1)}$.

\vspace{.1in}
\noindent
We use $\zeta_0$ to denote the point $\left ( \frac{1}{\sqrt{2}}, \frac{1}{\sqrt{2}}, ..., \frac{1}{\sqrt{2}}  \right ) $ in $\mathbb{R}^n$ and choose a sequence $\{t_k\}$ such that  $0< t_k < \frac{1}{4}$ and $\displaystyle \lim _{k \to \infty} t_k =0$.  Define the function $u(x)$ by  
\begin{equation}
	\label{defn:u}
	u(x) = \sum _{k=1}^{\infty} \epsilon_k r_k ^2 u_{t_k} \left ( \frac{x-R_k \zeta_0}{r_k} \right ) - \sum _{k=1}^{\infty} \epsilon_k r_k ^2 u_{t_k} \left ( \frac{x+R_k \zeta_0}{r_k} \right ) ,
\end{equation}
where the only conditions on $\epsilon_k$ for now are $\epsilon_k>0  $ and $\displaystyle \lim _{k \to \infty} \epsilon_k =0$,  we do not need to assign specific values to $\epsilon_k$ until later in the proof of Lemma \ref{lemma:u_Laplacian continuous_Hessian_bounded}.

\vspace{.05in}
\noindent
Condition (\ref{eqn:R_k_r_k}) ensures that the balls centered at the points $R_k\zeta_0$ with radii $r_k$ and  the balls centered at the points $- R_k\zeta_0$ with radii $r_k$ are all mutually disjoint. 

\begin{center}
	\includegraphics[scale=.3]{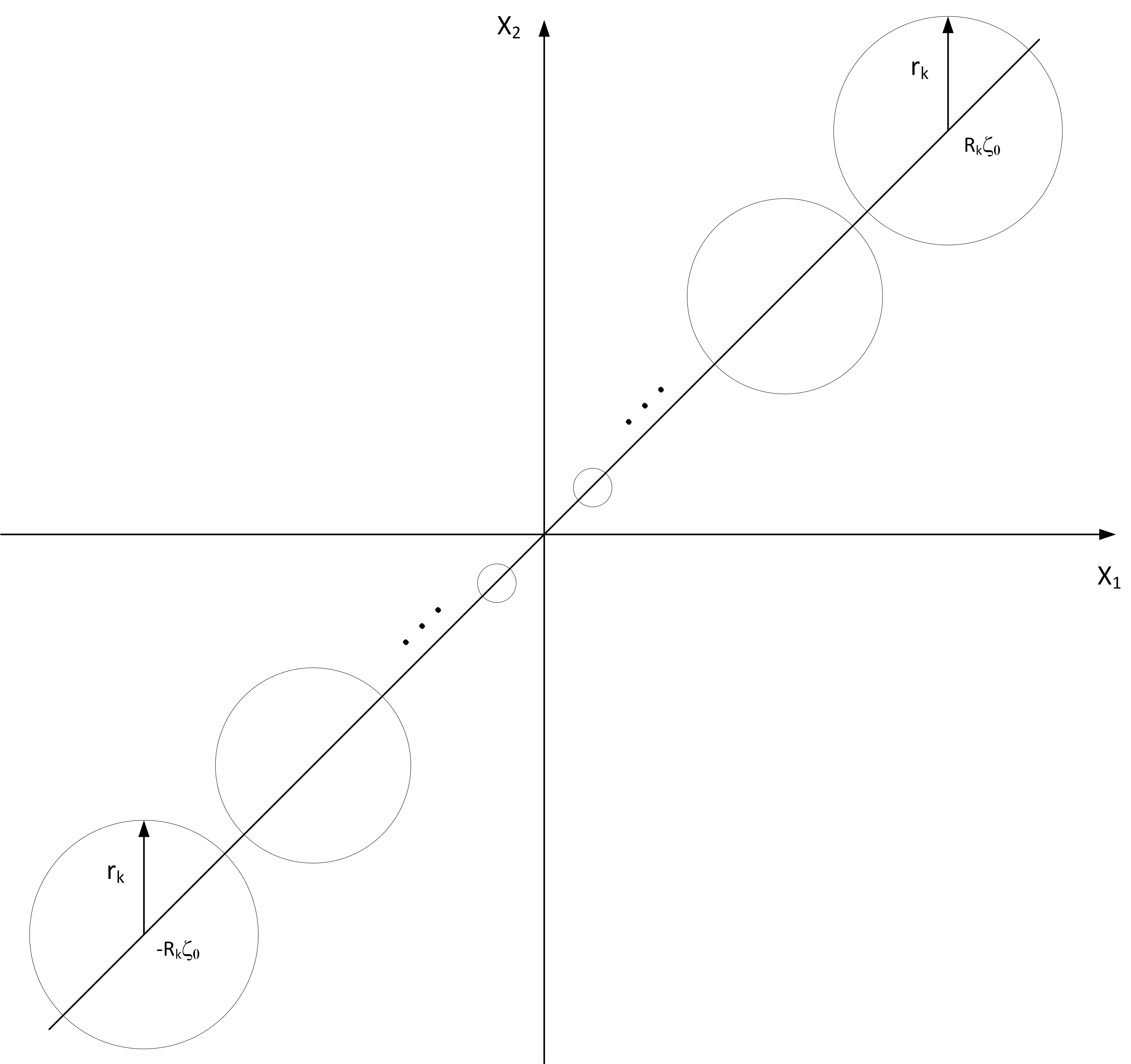}
\end{center}

\noindent
 For each $k \in \mathbb{N}$, let $B_k$ be the ball
centered at the point $R_k \zeta_0$ with radius $ \frac{2}{3}r_k$, and let $\tilde{B}_k $ be the ball
centered at the point $-R_k \zeta_0$ with radius $ \frac{2}{3}r_k$, then these $B_k$ and $\tilde{B}_k $ are also mutually disjoint.  By (\ref{eqn: condition_eta}) and (\ref{eqn: defn_u_t}), the support of each function $u_{t_k} \left ( \frac{x-R_k \zeta_0}{r_k} \right )$ is the ball $\{x \in \mathbb{R}^n: |x-R_k \zeta_0| \leq \frac{2}{3}r_k \}$, which is $B_k$, and the support of each function $u_{t_k} \left ( \frac{x+R_k \zeta_0}{r_k} \right )$ is the ball $\{x \in \mathbb{R}^n: |x+R_k \zeta_0| \leq \frac{2}{3}r_k \}$, which is $\tilde{B}_k$.
Therefore, although the definition of $u(x)$ appears to involve two infinite sums, it actually has only one non-zero term.   For any given $x \in \mathbb{R}^n$, if $x$ is outside of all the balls $B_k$ and $\tilde{B}_k$, then $$u(x)=0;$$ if $x \in B_k$ for some $k$, then $$u(x)=\epsilon_k r_k ^2 u_{t_k} \left ( \frac{x-R_k \zeta_0}{r_k} \right );$$  if $x \in \tilde{B}_k$ for some $k$, then $$u(x)= - \epsilon_k r_k ^2 u_{t_k} \left ( \frac{x + R_k \zeta_0}{r_k} \right ) .$$

\vspace{.05in}
\noindent
As $k \to \infty$, the radii of $B_k$ and $\tilde{B}_k$ go down to 0 and their centers move toward the origin, but none of the balls $B_k$ or $\tilde{B}_k$ contains the origin. In fact, for any $j=1,.., n$, the $x_j$ coordinate hyperplane does not intersect any of the balls $B_k$ or $\tilde{B}_k$.  To see this, let $$(x_1,...,x_{j-1}, 0, x_{j+1},...x_n)$$ be an arbitrary point on  the $x_j$ coordinate hyperplane. The distance from this point to the center of $B_k$, which is $R_k\zeta_0 = \left ( \frac{R_k}{\sqrt{2}}, ..., \frac{R_k}{\sqrt{2}}  \right ) $, is 

$$ \sqrt{ \left( \frac{R_k}{\sqrt{2}} -x_1 \right )^2+ \cdots + \left( \frac{R_k}{\sqrt{2}} -x_{j-1 }\right )^2+ \left( \frac{R_k}{\sqrt{2}} \right )^2 + \left( \frac{R_k}{\sqrt{2}} -x_{j+1} \right )^2 + \cdots + \left( \frac{R_k}{\sqrt{2}} -x_n \right )^2 } $$
$$ \geq   \frac{R_k}{\sqrt{2}}> \frac{r_k}{\sqrt{2}} > \frac{2}{3}r_k, \hspace{6in}
$$

\vspace{.05in}
\noindent
since $\frac{1}{\sqrt{2}} \approx 0.71$ and $ \frac{2}{3} \approx 0.67$.  Similarly, we can prove that the $x_j$ coordinate hyperplane does not intersect $\tilde{B}_k$ either. Thus $u = 0$ on all of the $n$ coordinate hyperplanes, and consequently $u(0)=0.$ By (\ref{eqn:u_t_bounded}) we know $u_{t_k}$ is uniformly bounded by a constant independent of $t_k$, which implies $\displaystyle \lim _{|x| \to 0} u(x)=0$.  Therefore $u$ is continuous at the origin, and hence continuous everywhere in $\mathbb{R}^n$.  Furthermore, the next lemma shows that $u$ is actually twice differentiable everywhere.

\vspace{.05in}
\begin{lemma}
	\label{lemma:u_differentiable}
	The function $u(x)$ as defined in (\ref{defn:u}) is twice differentiable everywhere in $\mathbb{R}^n$, and all its first and second order partial derivatives at the origin are equal to 0. 
\end{lemma}

\pf 
By definition $u(x)$ is $C^2$ for all $x \neq 0$, so we only need to show it is twice differentiable at the origin.
Because $u=0$ on all the coordinate hyperplanes, for any $i,j=1,...,n$,
\begin{equation*}
	\label{eqn:u_derivative_0}
	\frac{\partial u}{\partial x_j} (0)=0 \hspace{.4in} \text{and} \hspace{.4in} \frac{\partial^2 u}{\partial x_i \partial x_j} (0)=0.
\end{equation*}

\noindent
Thus
$$ \displaystyle \lim _{|x| \to 0}\frac { 	u (x) -u (0) - \displaystyle \sum _{i =1}^n \frac{\partial u}{\partial x_i} (0) x_i  -  
	\frac{1}{2}	\sum _{i,j=1}^n	\frac{\partial ^2 u}{\partial x_i \partial x_j}(0) x_i x_j  } {|x|^2}
= 
\lim _{|x| \to 0} \frac {u (x) } {|x|^2}. $$

\vspace{.05in}
\noindent
Recall that the balls $B_k$ and $ \tilde{B}_k$ are mutually disjoint, so for any $x \in \mathbb{R}^n$, if $x$ is outside of all the balls $B_k$ and $ \tilde{B}_k$, then  $$u(x)=0;$$ if $x \in B_k$ for some $k$, then $$u(x) = \epsilon_k r_k^2   u_{t_k} \left ( \frac{x-R_k \zeta_0}{r_k} \right );$$ if $x \in \tilde{B}_k$ for some $k$, then $$u(x) =  -\epsilon_k r_k^2   u_{t_k} \left ( \frac{x+R_k \zeta_0}{r_k} \right ) .$$

\vspace{.05in}
\noindent
By (\ref{eqn:u_t_bounded}) we have $|u_{t_k}| \leq C_{\eta, \varphi}$ which only depends on $\eta$ and $\varphi$. Thus

$$
\frac{|u(x) | }{|x|^2 }  \leq  \frac{\epsilon_k r_k^2 C_{\eta, \varphi} }{\left ( R_k - \frac{2}{3}r_k \right )^2} = \frac{\epsilon_k  C_{\eta, \varphi} }{\left ( \frac{R_k}{r_k} - \frac{2}{3}\right )^2 } <  \frac{\epsilon_k   C_{\eta, \varphi} }{\left ( 1 - \frac{2}{3} \right )^2 } \to 0  \hspace{.2in} \text{as} \,\, k \to \infty. 
$$

\vspace{.1in}
\noindent
Hence we know that $$ \lim _{|x| \to 0}\frac { 	u (x) -u (0) - \displaystyle \sum _{i =1}^n \frac{\partial u}{\partial x_i} (0) x_i  -  
	\frac{1}{2}	\sum _{i,j=1}^n	\frac{\partial ^2 u}{\partial x_i \partial x_j}(0) x_i x_j  } {|x|^2} =0, $$ 

\vspace{.1in}
\noindent
which means $u$ is twice differentiable at the origin. This completes the proof.

\stop

\vspace{.1in}
Next, we will show that $u$ has continuous Laplacian and bounded Hessian.

\vspace{.05in}
\begin{lemma}
	\label{lemma:u_Laplacian continuous_Hessian_bounded}
	The function $u(x)$ as defined in (\ref{defn:u}) has continuous Laplacian everywhere in $\mathbb{R}^n$, and all its second order partial derivatives are bounded.
\end{lemma}

\pf
Because $u=0$ on all the coordinate hyperplanes,
\begin{equation*}
	\label{eqn:u_derivative_0}
	\frac{\partial ^2 u}{\partial x_j^2} (0)=0.
\end{equation*}

\noindent
 For any given $x \in \mathbb{R}^n$, if $x$ is outside of all the balls $B_k$ and $ \tilde{B}_k$, then $$\frac{\partial^2 u}{\partial x_j^2}(x)=0;$$ if $x \in B_k$ for some $k$, then $$\frac{\partial ^2 u}{\partial x_j^2}(x) = \epsilon_k   \frac{\partial^2 u_{t_k}}{\partial x_j^2} \left ( \frac{x-R_k \zeta_0}{r_k} \right ) ;$$ if $x \in \tilde{B}_k$ for some $k$, then $$\frac{\partial ^2 u}{\partial x_j^2}(x) = - \epsilon_k   \frac{\partial^2 u_{t_k}}{\partial x_j^2} \left ( \frac{x+R_k \zeta_0}{r_k} \right ).$$ 
 
 \noindent
 By (\ref{eqn:u_t_derivative_bounded_second_diagonal}) in Lemma \ref{lemma:u_t_derivative_bounded}, 
 \begin{equation}
 	\label{eqn: ut_xj^2_positive}
    \epsilon _k \left \vert  \frac{\partial^2 u_{t_k}}{\partial x_j^2} \left ( \frac{x-R_k \zeta_0}{r_k} \right )   \right \vert  \leq \epsilon _k   C_{\eta, \varphi} \to 0 \hspace{.2in} \text{as} \,\, k \to \infty ,
\end{equation}
     and 
     \begin{equation}
     	\label{eqn: ut_xj^2_negtive}
        \epsilon _k \left \vert  \frac{\partial^2 u_{t_k}}{\partial x_j^2} \left ( \frac{x+R_k \zeta_0}{r_k} \right )   \right \vert  \leq \epsilon _k   C_{\eta, \varphi} \to 0 \hspace{.2in} \text{as} \,\, k \to \infty. 
    \end{equation}
           Thus $\displaystyle \lim_{|x| \to 0} \frac{\partial ^2 u}{\partial x_j^2}(x) = 0 $, which implies that $\displaystyle \frac{\partial^2 u}{\partial x_j^2}(x)$ is continuous at 0.  Since it is also continuous for all $x \neq 0$, it is continuous everywhere.
This proves that $\displaystyle \Delta u = \sum _{j=1}^n \frac{\partial^2 u}{\partial x_j^2}$ is continuous everywhere in $\mathbb{R}^n$.

\vspace{.1in}
\noindent
Next we will show that all the second derivatives of $u$ are bounded.  By (\ref{eqn: ut_xj^2_positive}) and (\ref{eqn: ut_xj^2_negtive}) we can conclude that $\frac{\partial^2 u}{\partial x_j^2} $ is bounded for all $j=1,... n$, thus we only need to show that the non-diagonal entries of $\Hess u$ are also bounded.

 \vspace{.05in}
\noindent
For any $x \in \mathbb{R}^n$ and $j=3,...,n$,  there are three possible cases for $\frac{\partial ^2 u}{\partial x_1 \partial x_j} (x)$:

\noindent
if $x$ is outside of all the balls $B_k$ and $ \tilde{B}_k$, then $$\frac{\partial^2 u}{\partial x_1 \partial x_j}(x)=0;$$ if $x \in B_k$ for some $k$, then $$\frac{\partial ^2 u}{\partial x_1 \partial x_j}(x) = \epsilon_k   \frac{\partial^2 u_{t_k}}{\partial x_1 \partial x_j} \left ( \frac{x-R_k \zeta_0}{r_k} \right );$$ if $x \in \tilde{B}_k$ for some $k$, then $$\frac{\partial ^2 u}{\partial x_1 \partial x_j}(x) = - \epsilon_k   \frac{\partial^2 u_{t_k}}{\partial x_1 \partial x_j} \left ( \frac{x+R_k \zeta_0}{r_k} \right ).$$

\noindent
By (\ref{eqn:u_t_derivative_bounded_x_1_x_j}) in Lemma \ref{lemma:u_t_derivative_bounded_non-diagonal},  $\frac{\partial ^2 u_{t_k}}{\partial x_1 \partial x_j} (x)$ is bounded by a constant independent of $t_k$.  Then since $\epsilon_k \to 0 $, we know that $\frac{\partial ^2 u}{\partial x_1 \partial x_j} (x) \to 0$ as $x \to 0$.  In particular, it is bounded by a constant depending only on $\eta$ and $\varphi$.  Similarly, we can prove that $\frac{\partial ^2 u}{\partial x_2 \partial x_j} (x)$ and $\frac{\partial ^2 u}{\partial x_i \partial x_j} (x)$, where  $i,j=3,...,n$ and $i \neq j$, also approach $0$ as $x \to 0$.  In particular, they are also bounded.

\vspace{.1in}
\noindent
The last derivative we need to look at is $\frac{\partial ^2 u}{\partial x_1 \partial x_2} (x)$. The proof of its boundedness will be different from that for the other derivatives.
\allowdisplaybreaks
\begin{eqnarray}
	\label{eqn:2nd_derivatives_u_t_x_1_x_2}
& &	\frac{\partial ^2 u_t}{\partial x_1 \partial x_2} (x) \\
 & = &  \eta  (|x|)    |x|^{2t} \varphi (-\ln  |x|^2 )  + \eta '' (|x|)   x_1^2 x_2 ^2 |x|^{2t-2} \varphi  (-\ln  |x|^2   )  \nonumber \\
	& + &   \eta ' (|x|) \left ( x_1^2 + x_2^2 \right )  |x|^{2t-1} \varphi  (-\ln  |x|^2   ) + (4t-1)  \eta ' (|x|)   x_1^2 x_2^2|x|^{2t-3} \varphi  (-\ln  |x|^2   ) \nonumber  \\
	& + &    2t(2t-2) \eta  (|x|)   x_1 ^2 x_2^2 |x|^{2t-4} \varphi (-\ln  |x|^2  ) + 2t \eta  (|x|) \left ( x_1^2 + x_2^2 \right )    |x|^{2t-2} \varphi (-\ln  |x|^2  )  \nonumber \\
&	- &  4 \eta ' (|x|)   x_1^2 x_2^2 |x|^{2t-3} \varphi '  (-\ln  |x|^2  ) - 2 \eta  (|x|) \left ( x_1^2 + x_2^2 \right )   |x|^{2t-2} \varphi ' (-\ln  |x|^2  )     \nonumber \\
& + &  (4-8t) \eta  (|x|)   x_1^2 x_2^2 |x|^{2t-4} \varphi '  (-\ln |x|^2  )    +  4 \eta  (|x|)   x_1^2 x_2^2  |x|^{2t-4} \varphi '' (-\ln |x|^2  ). \nonumber
\end{eqnarray}

\vspace{.1in}
\noindent
The first term, $ \eta  (|x|)    |x|^{2t} \varphi (-\ln  |x|^2 ) $, is bounded by $C |x|^{2t} \left \vert \varphi (-\ln  |x|^2 ) \right \vert .$ All the other terms in (\ref{eqn:2nd_derivatives_u_t_x_1_x_2}) are bounded by a constant independent of $t$, by the same estimates as those in the proof of Lemma \ref{lemma:u_t_derivative_bounded_non-diagonal}.  Therefore, we have 
\begin{equation}
	\label{eqn:ut_x1_x2_bound}
    \left \vert \frac{\partial ^2 u_t}{\partial x_1 \partial x_2} (x)	\right \vert  \leq  C \sup_{|x| \leq \frac{2}{3} } \left ( |x|^{2t} \left \vert \varphi (-\ln  |x|^2 ) \right \vert \right ) +  C_{\eta, \varphi},
\end{equation} 
where $C$ and $ C_{\eta, \varphi}$ are constants depending on $\eta$ and $\varphi$, and we have used the fact that by definition $\eta(x) =0$ when $|x| > \frac{2}{3}$.

\vspace{.1in}
\noindent
For any $x \in \mathbb{R}^n$,  there are three possible cases for $\frac{\partial ^2 u}{\partial x_1 \partial x_2} (x)$:

\vspace{.05in}
\noindent
if $x$ is outside of all the balls $B_k$ and $ \tilde{B}_k$, then $$\frac{\partial^2 u}{\partial x_1 \partial x_2}(x)=0;$$ if $x \in B_k$ for some $k$, then $$\frac{\partial ^2 u}{\partial x_1 \partial x_2}(x) = \epsilon_k   \frac{\partial^2 u_{t_k}}{\partial x_1 \partial x_2} \left ( \frac{x-R_k \zeta_0}{r_k} \right );$$ if $x \in \tilde{B}_k$ for some $k$, then $$\frac{\partial ^2 u}{\partial x_1 \partial x_2}(x) = - \epsilon_k   \frac{\partial^2 u_{t_k}}{\partial x_1 \partial x_2} \left ( \frac{x+R_k \zeta_0}{r_k} \right ).$$ 

\vspace{.05in}
\noindent
Thus in any case, by (\ref{eqn:ut_x1_x2_bound}) we have 

$$  \left \vert \frac{\partial ^2 u}{\partial x_1 \partial x_2}(x)  \right \vert \leq \ \epsilon_k \left ( C \sup_{|x| \leq \frac{2}{3} } \left (|x|^{2t_k} \left \vert \varphi (-\ln  |x|^2 ) \right \vert \right ) +  C_{\eta, \varphi} \right ) .$$

\vspace{.05in}
\noindent
Denote 
\begin{equation}
	\label{defn:M_k}
	 M_k = \sup _{|x| \leq \frac{2}{3}} \left ( |x|^{2t_k}  \varphi (-\ln  |x|^2 )  \right ) ,
\end{equation}
  then  $$  \left \vert \frac{\partial ^2 u}{\partial x_1 \partial x_2}(x)  \right \vert \leq \ C\epsilon_k M_k + \epsilon_k  C_{\eta, \varphi} .$$
Since $t_k \to 0$, by Lemma \ref{lem: x-time-varphi-infinity} we know $\displaystyle \sup_{|x| \leq \frac{1}{2}} \left ( |x|^{2t_k}  \varphi (-\ln  |x|^2 ) \right ) \to \infty$, and consequently  $M_k \to \infty$. Now we let $$\epsilon_k = \frac{1}{M_k} ,  $$ then 
 $$  \left \vert \frac{\partial ^2 u}{\partial x_1 \partial x_2}(x)  \right \vert \leq \ C +  C_{\eta, \varphi} ,$$ thus $\frac{\partial ^2 u}{\partial x_1 \partial x_2} $ is bounded.  This completes the proof of Lemma \ref{lemma:u_Laplacian continuous_Hessian_bounded}.

\stop

\vspace{.1in}

To finish the proof of Theorem \ref{thm:bounded_2nd_order} we only need to prove the following lemma.

\vspace{.05in}
\noindent

\begin{lemma}	
	\label{lemm: u_not_C2}
	The function $u$ defined by (\ref{defn:u}) is not in  $C^2(\mathbb{R}^n)$.
\end{lemma}

\vspace{.05in}
\noindent
\pf
we will prove that $\frac{\partial ^2 u}{\partial x_1 \partial x_2}$ is not continuous at the origin.  Recall that in Lemma \ref{lemma:u_differentiable} we showed $\frac{\partial ^2 u}{\partial x_1 \partial x_2}(0)=0$, so we need to show that $\displaystyle \lim_{x \to 0} \frac{\partial ^2 u}{\partial x_1 \partial x_2}(x)$  is not 0.
\vspace{.05in}
\noindent
For each $k$, by Lemma \ref{lem: x-time-varphi-infinity} we may choose a point $y^{(k)}$, where $|y^{(k)}|\leq \frac{1}{2}$, such that 
\begin{equation}
	\label{eqn:choose_y(k)}
	|y^{(k)}|^{2t_k}  \varphi (-\ln  |y^{(k)}|^2 )  = \frac{1}{2}M_k.
\end{equation}

\vspace{.05in}
\noindent
Let $x^{(k)} \in \mathbb{R}^n$ be the point satisfying
\begin{equation}
	\label{eqn:defn_x_(k)}
	\frac{x^{(k)} - R_k\zeta _0}{r_k} =  y^{(k)}, \nonumber
\end{equation}
and let $\tilde{x}^{(k)} \in \mathbb{R}^n$ be the point satisfying
\begin{equation}
	\label{eqn:defn_x_(k)}
	\frac{\tilde{x}^{(k)} + R_k\zeta _0}{r_k} =y^{(k)}. \nonumber
\end{equation}

\vspace{.05in}
\noindent
Then $ x^{(k)} \in B_k$, $ \tilde{x}^{(k)} \in \tilde{B}_k, $ and $\displaystyle \lim _{k\to \infty}  x^{(k)} = \lim _{k\to \infty} \tilde{x}^{(k)} =0.  $ By the definition of $u$,

\begin{eqnarray*}
\frac{\partial ^2 u}{\partial x_1 \partial x_2}( x^{(k)}) & = & \epsilon_k   \frac{\partial^2 u_{t_k}}{\partial x_1 \partial x_2} \left ( \frac{ x^{(k)}-R_k \zeta_0}{r_k} \right )	\\
& = & \epsilon_k   \frac{\partial^2 u_{t_k}}{\partial x_1 \partial x_2} \left (  y^{(k)}\right ),
\end{eqnarray*}
and
\begin{eqnarray*}
	\frac{\partial ^2 u}{\partial x_1 \partial x_2}( \tilde{x}^{(k)}) & =  & - \epsilon_k   \frac{\partial^2 u_{t_k}}{\partial x_1 \partial x_2} \left ( \frac{ \tilde{x}^{(k)} + R_k \zeta_0}{r_k} \right )	\\
	& = & - \epsilon_k   \frac{\partial^2 u_{t_k}}{\partial x_1 \partial x_2} \left (  y^{(k)}\right ).
\end{eqnarray*}

\vspace{.05in}
\noindent
By (\ref{eqn:2nd_derivatives_u_t_x_1_x_2}) and (\ref{eqn:choose_y(k)}), we have 
\begin{eqnarray*}
	 \frac{\partial^2 u_{t_k}}{\partial x_1 \partial x_2} \left (  y^{(k)} \right ) & = &   \eta  (|y^{(k)}|)    |y^{(k)}|^{2t_k} \varphi (-\ln  |y^{(k)}|^2 )  + h(y^{(k)}) \\
	 & = & \frac{1}{2}   M_k  + h(y^{(k)}),
\end{eqnarray*}
where the $h$ function is the sum of the second to last terms in (\ref{eqn:2nd_derivatives_u_t_x_1_x_2}), and we have used the fact that $  \eta  (|y^{(k)}|) =1$ since $|y^{(k)}| \leq \frac{1}{2}$.  

\vspace{.05in}
\noindent
Thus we know 
\begin{equation*}
	\label{eqn:derivative_u_x1_x2_first_sequence}
\frac{\partial ^2 u}{\partial x_1 \partial x_2}( x^{(k)}) = \epsilon_k \left ( \frac{1}{2}   M_k  + h(y^{(k)}) \right ) = 	\frac{1}{2} + \epsilon_k h(y^{(k)}),
\end{equation*}
and 
\begin{equation*}
	\label{eqn:derivative_u_x1_x2_second_sequence}
	\frac{\partial ^2 u}{\partial x_1 \partial x_2}( \tilde{x}^{(k)}) = - \epsilon_k \left ( \frac{1}{2}   M_k  + h(y^{(k)}) \right ) = -	\frac{1}{2} - \epsilon_k h(y^{(k)}).
\end{equation*}

\vspace{.05in}
\noindent
As shown in the proof of Lemma \ref{lemma:u_Laplacian continuous_Hessian_bounded}, $h$ is bounded by a constant $C_{\eta, \varphi}$ that depends only on $\eta$ and $\varphi$.  Thus $\displaystyle \lim_{k \to \infty} \epsilon_k h(y^{(k)}) =0,$ and consequently   $$ \lim_{k \to \infty} \frac{\partial ^2 u}{\partial x_1 \partial x_2}( x^{(k)}) = \frac{1}{2} \hspace{.2in} \text{and}  \hspace{.2in} \lim_{k \to \infty} \frac{\partial ^2 u}{\partial x_1 \partial x_2}( \tilde{x}^{(k)}) = - \frac{1}{2}.   $$
This implies that $\frac{\partial ^2 u}{\partial x_1 \partial x_2} $ is not continuous at $0$, so $u$ is not $C^2$.

\stop


\vspace{.2in}
The idea for constructing higher order examples for Theorem \ref{thm: bounded_kth_order} is the same as that for Theorem \ref{thm:bounded_2nd_order}.  As in the higher order case in \cite{Pan-Yan_unbounded}, to simplify the calculations we use a complex variable for the first two components of $x$: for any $x=(x_1, x_2, x_3,...,x_n) \in \mathbb{R}^n$, denote
 \begin{equation*}
 	\label{eqn:define_z}
 	z=x_1+ix_2 \,\, \,\, \text{and} \,\, \,\, \bar{z}=x_1-ix_2.
 \end{equation*} 
Then
$$ x_1^2+x_2^2=z\bar{z}, \hspace{.4in} |x|^2=z\bar{z}+\sum_{j=3}^n x_j ^2, \hspace{.4in} \text{and} \hspace{.4in} \frac{\partial ^2 }{\partial x_1^2} + 	\frac{\partial ^2 }{\partial x_2^2} =4	\frac{\partial ^2 }{\partial \bar{z} \partial z},$$ and consequently
\begin{eqnarray*}
	\frac{\partial |x|}{ \partial z} = \frac{\bar{z}}{2|x|}, \hspace{.4in} \frac{\partial |x|}{ \partial \bar{z}} = \frac{z}{2|x|}, \hspace{.4in} \text{and} \hspace{.4in} \frac{\partial |x|}{ \partial x_j} = \frac{x_j}{|x|} \,\, \,\, (\text{when} \,\, j \geq 3).
\end{eqnarray*}  
 
 
\vspace{.05in}
 \noindent
 Recall that in the higher order case $\varphi(s)$ needs to be $(k+2)$-times differentiable and satisfy
 \begin{equation}
 	\label{eqn:higher_condition_phi}
 		\lim _{s \to \infty} \varphi  (s) = \infty,  \hspace{.2in} \lim _{s \to \infty} \varphi ' (s) = \cdots = \lim _{s \to \infty} \varphi ^{(k+2)}(s) = 0.
 \end{equation}
 
 
\vspace{.05in}
 \noindent
The function $u_t: \mathbb{R}^n \to \mathbb{R} \,\, (n \geq 2)$ is now defined by
\begin{equation}
	\label{eqn:higher_order_defn_u_t}
	u_t(x) = 	\left\{
	\begin{array}{l l l}      
		0 &   \hspace{.1in} x=0, \\
		\noalign{\medskip}
		\eta(|x|)z^{k+2}|x|^{2t}\varphi  (-\ln  |x|^2  ) & \hspace{.1in} 0 < |x|< 1, \\
		\noalign{\medskip}
		0 &	\hspace{.1in} |x| \geq 1,
	\end{array}\right.
\end{equation}
where $\eta$ is the same as in (\ref{eqn: condition_eta}) and $\varphi$ satisfies (\ref{eqn:higher_condition_phi}).  Then we define $u$ in the same way as in (\ref{defn:u}): 
$$u (x)= \sum _{k=1}^{\infty} \epsilon_k r_k ^2 u_{t_k} \left ( \frac{x-R_k \zeta_0}{r_k} \right ) - \sum _{k=1}^{\infty} \epsilon_k r_k^2   u_{t_k} \left ( \frac{x+R_k \zeta_0}{r_k} \right ) . $$

\vspace{.1in}
\noindent
By similar but lengthier calculations we can show that  the complex-valued function $u$ is $(k+2)$-times differentiable, $\Delta u$ is $C^k$ throughout $\mathbb{R}^n$, $D^{k+2}u$ is bounded, but $u$ is not in $C^{k+2}$. The real and imaginary parts of $u$ are two real-valued functions that are $(k+2)$-times differentiable at 0, their Laplacian are $C^k$ throughout $\mathbb{R}^n$, and their $(k+2)$-th partial derivatives are bounded, but at least one of their $(k+2)$-th derivatives is not continuous at the origin.  Therefore, we have found a real-valued function that satisfies all the conditions in Theorem \ref{thm: bounded_kth_order}.

\vspace{.2in}

\section{Examples for the Monge-Amp\`ere Equation}
\label{section:monge-ampere}

\vspace{.1in}
Moving from Theorem \ref{thm:bounded_2nd_order} to Theorem \ref{thm: monge-ampere}, the only difference in the conclusion is that we need to show $\det D^2u$ is continuous instead of $\Delta u$ being continuous.  We will define $u$ in the same way as in the proof of Theorem \ref{thm:bounded_2nd_order}, but for technical reasons that will be explained later, the proof only works when the dimension of $\mathbb{R}^n$ is $n \geq 3$, and we need to assume the following additional conditions on $\varphi$:
\begin{equation}
	\label{eqn:condition_phi_monge_ampere}
	\lim_{s \to \infty} \varphi ^{n-1} (s) \varphi '(s) =0   \hspace{.2in} \text{and} \hspace{.2in} 		\lim_{s \to \infty} \varphi ^{n-1} (s) \varphi ''(s) =0  .
	\end{equation}

Let $u$ be the function defined in (\ref{defn:u}).  By definition $\det D^2u(x)$ is continuous for all $x \neq 0$.  As proved in Section \ref{section:2nd_proof}, all the partial derivatives of $u$ is equal to 0 at the origin, so $\det D^2 u(0)=0$. Therefore we only need to show $\displaystyle \lim _{|x| \to 0} \det D^2 u(x)=0 $, which can be proved by an argument similar to that for the continuity of $\Delta u$ in the proof of Lemma \ref{lemma:u_Laplacian continuous_Hessian_bounded}.  As in the previous section, for each $k \in \mathbb{N}$, let $B_k$ be the ball
centered at the point $R_k \zeta_0$ with radius $ \frac{2}{3}r_k$, and let $\tilde{B}_k $ be the ball
centered at the point $-R_k \zeta_0$ with radius $ \frac{2}{3}r_k$.  For any $x \in \mathbb{R}^n$,  there are three possible cases for $D^2 u (x)$:

\vspace{.05in}
\noindent
if $x$ is outside of all the balls $B_k$ and $ \tilde{B}_k$, then $$D^2u(x)=0;$$ if $x \in B_k$ for some $k$, then $$ D^2 u(x) = \epsilon_k   D^2 u_{t_k} \left ( \frac{x-R_k \zeta_0}{r_k} \right );$$ if $x \in \tilde{B}_k$ for some $k$, then $$ D^2 u(x) = - \epsilon_k   D^2 u_{t_k} \left ( \frac{x+R_k \zeta_0}{r_k} \right ).$$ 

\vspace{.05in}
\noindent
Since $\epsilon_k \to 0$, to show $\displaystyle \lim _{|x| \to 0} \det D^2 u(x)=0, $ we only need to show $\det D^2 u_t$ is uniformly bounded by a constant independent of $t$. 

\vspace{.05in}
In the proofs of Lemma \ref{lemma:u_t_derivative_bounded}, Lemma \ref{lemma:u_t_derivative_bounded_non-diagonal}, and Lemma \ref{lemma:u_Laplacian continuous_Hessian_bounded}, all the second order partial derivatives of $u_t$ were estimated term-by-term, and we have the following estimates. (In the subsequent discussions in this section we will use $C$ to denote a constant depending only on $\eta$ and $\varphi$ )

\begin{eqnarray}
	\label{eqn:second_order_bound_x1_x_2}
 \left \vert	\frac{\partial ^2 u_t}{\partial x_1 \partial x_2} (x) \right \vert 
& \leq & \sup _{|x| \leq \frac{2}{3}} \biggl(  C|x|^{2t}|\varphi(-\ln|x|^2) | + C t|x|^{2t} |\varphi (-\ln|x|^2) | +  C |x| |\varphi (-\ln|x|^2) | \nonumber \\
& +  &  C|\varphi '(-\ln|x|^2) | +  C|\varphi ''(-\ln|x|^2) | \biggr ) ,  
\end{eqnarray}
\begin{eqnarray}
	\label{eqn:second_order_bound_xi^2}
		\left \vert	\frac{\partial ^2 u_t}{\partial x_i ^2} (x) \right \vert  & \leq & \sup _{|x| \leq \frac{2}{3}} \biggl(  C t|x|^{2t} |\varphi (-\ln|x|^2) | +  C |x| |\varphi (-\ln|x|^2) | \nonumber \\
		& +  &  C|\varphi '(-\ln|x|^2) | +  C|\varphi ''(-\ln|x|^2) |  \biggr ),     \hspace{1in}  i=1,2,..., n, 
		\end{eqnarray}
		\begin{eqnarray}
			\label{eqn:second_order_bound_x1_x_i}
		\left \vert	\frac{\partial ^2 u_t}{\partial x_1 x_j} (x) \right \vert 
	& \leq & \sup _{|x| \leq \frac{2}{3}} \biggl( C t|x|^{2t} |\varphi (-\ln|x|^2) | +  C |x| |\varphi (-\ln|x|^2) | \nonumber \\
	& +  &  C|\varphi '(-\ln|x|^2) | +  C|\varphi ''(-\ln|x|^2) |  \biggr ),   \hspace{1in}  j=3,..., n, 
\end{eqnarray}
\begin{eqnarray}
	\label{eqn:second_order_bound_xi_x_j}
	\left \vert	\frac{\partial ^2 u_t}{\partial x_i x_j} (x) \right \vert 
	& \leq & \sup _{|x| \leq \frac{2}{3}} \biggl(  C t|x|^{2t} |\varphi (-\ln|x|^2) | +  C |x| |\varphi (-\ln|x|^2) | \nonumber \\
	& +  &  C|\varphi '(-\ln|x|^2) | +  C|\varphi ''(-\ln|x|^2) |  \biggr ),   \hspace{1in}  i,j=3,..., n, \, \text{and} \,\, i \neq j. 
\end{eqnarray}

\noindent
Thus we only need to consider $|x| \leq \frac{2}{3}$.  For simplicity we first assume $n=3$.  Then $$\det D^2u_t = \det	\left [
\begin{array} {ccc}
	\frac{\partial ^2 u_t}{\partial x_1 ^2}(x)  & \frac{\partial ^2 u_t}{\partial x_1 \partial x_2} (x) &   \frac{\partial ^2 u_t}{\partial x_1 \partial x_3} (x) \\
	\noalign{\bigskip}
	\frac{\partial ^2 u_t}{\partial x_2 \partial x_1 } (x) & \frac{\partial ^2 u_t}{\partial x_2^2} (x) 	&  \frac{\partial ^2 u_t}{\partial x_2 \partial x_3} (x) \\
	\noalign{\bigskip}
	\frac{\partial ^2 u_t}{\partial x_3 \partial x_1 } (x) & \frac{\partial ^2 u_t}{\partial x_3 \partial x_2} (x) 	& \frac{\partial ^2 u_t}{\partial  x_3 ^2} (x)
\end{array}
\right ]. $$

\vspace{.05in}
\noindent
Applying (\ref{eqn:second_order_bound_x1_x_2}) to  (\ref{eqn:second_order_bound_xi_x_j}) in the expansion of this determinant, 
 and using the fact that $|x|<1$ and $0<t<1$ to simplify some terms, we have the estimate
 \vspace{.05in}
\allowdisplaybreaks
\begin{eqnarray}
	\label{eqn:bound_det_Hess}
 \left \vert \det D^2u_t \right \vert
	& \leq  &  C t|x|^{6t}|\varphi ^3 (-\ln|x|^2) |  +    C|x||\varphi ^3 (-\ln|x|^2) | + C|\varphi ^2 (-\ln|x|^2) | |\varphi ' (-\ln|x|^2) | \nonumber \\
	& + &  C|\varphi ^2 (-\ln|x|^2) | |\varphi '' (-\ln|x|^2) | +  C|\varphi  (-\ln|x|^2) | |\varphi ' (-\ln|x|^2) |^2   \nonumber \\
	& + &  C|\varphi  (-\ln|x|^2) | |\varphi '' (-\ln|x|^2) |^2 + C|\varphi ' (-\ln|x|^2) |^3 +   C|\varphi '' (-\ln|x|^2) |^3  \\
	& + &  C |\varphi ' (-\ln|x|^2) |^2|\varphi '' (-\ln|x|^2) | + C |\varphi ' (-\ln|x|^2) ||\varphi '' (-\ln|x|^2) |^2  . \nonumber
\end{eqnarray}

\vspace{.05in}
 \noindent
 By the conditions on $\varphi$ in (\ref{eqn: condition_phi}) and (\ref{eqn:condition_phi_monge_ampere}), and by Lemma \ref{lemma:x-times-varphi}, the third to last terms in (\ref{eqn:bound_det_Hess}) are all bounded by a constant independent of $t$.  As for the first two terms,
 denote $\Phi (s) = \varphi ^3 (s)$. Then $ \displaystyle \lim _{s \to \infty} \Phi  (s) = \infty$, and (\ref{eqn:condition_phi_monge_ampere}) implies that 
 $$ \lim _{s \to \infty} \Phi ' (s) =0 .  $$
 Thus the proofs of Lemma \ref{lemma:x-times-varphi} and Lemma \ref{lemma:t_times_x^t_time_varphi} still hold with $\varphi$ replaced by $\Phi$.  By Lemma \ref{lemma:t_times_x^t_time_varphi} we know that 
 $$t|x|^{6t}|\varphi ^3 (-\ln|x|^2) | = t|x|^{6t}|\Phi  (-\ln|x|^2) | $$ is bounded by a constant depending only on $\varphi$. By Lemma \ref{lemma:x-times-varphi} we know that  $$|x||\varphi ^3 (-\ln|x|^2) | = |x||\Phi (-\ln|x|^2) |  $$ has a removable discontinuity at 0, and hence it is bounded on the closed set $|x| \leq \frac{2}{3}$. Therefore, we have proved that $| \det D^2 u_t |$ is uniformly bounded by a constant independent of $t$, which then implies that $\det D^2 u$ is continuous.  This completes the proof of Theorem \ref{thm: monge-ampere} when $n=3$.  For $n \geq 4$, we denote $\Phi(s) = \varphi^n (s)$ and the estimates are essentially the same.
 
 \stop
 
 \vspace{.1in}
\noindent
\textbf{Remark 3.1.} The above proof would not hold in $\mathbb{R}^2$, because if $n=2$, then 
 $$\det D^2u_t = \det	\left [
 \begin{array} {cc}
 	\frac{\partial ^2 u_t}{\partial x_1 ^2}  & \frac{\partial ^2 u_t}{\partial x_1 \partial x_2}  \\
 	\noalign{\bigskip}
 	\frac{\partial ^2 u_t}{\partial x_2 \partial x_1 }  & \frac{\partial ^2 u_t}{\partial x_2^2} 
 \end{array}
 \right ] = 	\frac{\partial ^2 u_t}{\partial x_1 ^2} \cdot \frac{\partial ^2 u_t}{\partial x_2^2}  -  \left ( \frac{\partial ^2 u_t}{\partial x_1 \partial x_2}  \right )^2.$$
 \noindent
 While $\frac{\partial ^2 u_t}{\partial x_1 ^2}$ and $\frac{\partial ^2 u_t}{\partial x_2 ^2}$ are bounded independent of $t$, one of the terms in $\frac{\partial ^2 u_t}{\partial x_1 \partial x_2} $ is 
  $|x|^{2t}|\varphi (-\ln|x|^2)|$.  By Lemma \ref{lem: x-time-varphi-infinity} this term approaches $\infty$ as $|x| \to 0$ and $t \to 0$, thus $\det D^2 u_t$ is not uniformly bounded, and we cannot conclude that $\det D^2 u$ is continuous at the origin.

\vspace{.1in}
\noindent
\textbf{Remark 3.2.} As proved in Section \ref{section:2nd_proof}, the function $u$ in the above construction has bounded Hessian.  In \cite{Pan-Yan_unbounded}, we constructed another function $u$ that is twice differentiable everywhere with continuous Laplacian and unbounded Hessian. It is defined by 
 $$ 	u(x) = \sum _{k=1}^{\infty} \epsilon_k r_k ^2 u_{t_k} \left ( \frac{x-R_k \zeta_0}{r_k} \right ), $$
where $R_k$, $r_k$, and $\zeta_0$ are the same as in this paper, but we chose a different sequence of $\epsilon_k$.  For that $u$, by similar argument we can prove $\det D^2u$ is also continuous. Therefore, there is an ample supply of functions, some with bounded Hessian and some with unbounded Hessian, that are twice differentiable everywhere with continuous $\det D^2u$, but the functions are not $C^2$.

\vspace{.2in}

\section{The Geometric Application}
\label{section:geometric}

\vspace{.1in}

To prove Corollary \ref{cor:hypersurface}, let $u$ be a function constructed in Theorem \ref{thm:bounded_2nd_order}, and let $\Sigma$ be the hypersurface in $\mathbb{R}^{n+1}$ that is the graph of $u(x)$. The mean curvature of $\Sigma$ is 
$$ H=  \diver \left ( \frac{\triangledown u}{\sqrt{1+|\triangledown u|^2}} \right ) = \frac{\Delta u}{\sqrt{1+|\triangledown u|^2}} - \sum_{i,j=1}^n \frac{\frac{\partial u}{\partial x_i} \cdot \frac{\partial u}{\partial x_j} \cdot \frac{\partial^2 u}{\partial x_i \partial x_j}   }{\sqrt{1+|\triangledown u|^2}^3} .  $$ Since $\Delta u$ and $\triangledown u$ are both continuous, $ \frac{\Delta u}{\sqrt{1+|\triangledown u|^2}} $ is continuous.  As shown in the proof of Theorem \ref{thm:bounded_2nd_order}, $ \frac{\partial u}{\partial x_i}$ and $ \frac{\partial u}{\partial x_j}$ both approach 0 as $|x| \to 0$.  Then since $\frac{\partial^2 u}{\partial x_i \partial x_j}$ is bounded, $\frac{\partial u}{\partial x_i} \frac{\partial u}{\partial x_j}  \frac{\partial^2 u}{\partial x_i \partial x_j}   $ approaches $0$ as $|x| \to 0$. As a result $H$ is continuous at the origin, hence continuous everywhere.  On the other hand, since $u$ is not $C^2$, the hypersurface $\Sigma$ is not $C^2$.  The second fundamental form of $\Sigma$ is given by the matrix 
$$ \frac{\Hess u}{\sqrt{1+|\triangledown u|^2}}. $$ 
It is bounded but not continuous since $\Hess u$ is bounded and discontinuous.  This proves Corollary \ref{cor:hypersurface}.

\vspace{.2in}

\section{The Radially Symmetric Case}
\label{section:radial}

\vspace{.1in}
In this section we prove Proposition \ref{prop:radially_symmetric}.  Let $\omega(x) $ be a twice differentiable, radially symmetric function. First of all, without loss of generality we can assume $\omega(0)=0$,  therefore we can express $\omega$ as
\begin{equation*}
	\omega (x) =	\left\{
	\begin{array}{l l }      
		\psi(|x|) & \hspace{.4in}  x \neq (0,...,0) , \\
		\noalign{\smallskip}
		0 & \hspace{.4in}   x = (0,...,0) ,
	\end{array}\right.
\end{equation*}
where $\psi(s): [0, \infty) \to \mathbb{R}$ is twice differentiable.  Next, we make the observation that $\psi '(0)=0$.  To see this, note that 
\begin{eqnarray*}
	\frac{\partial \omega}{\partial x_1}(0) & = & \lim _{h \to 0} \frac{\psi(|h|)}{h} \\
	& = & \lim _{h \to 0} \frac{\psi(|h|)}{|h|} \frac{|h|}{h}\\
	& = & \pm \psi'(0),
\end{eqnarray*}
so to ensure the existence of $ \frac{\partial \omega}{\partial x_1}(0) $ we must have $\psi '(0)=0$.  Consequently, Proposition \ref{prop:radially_symmetric} is a direct corollary of the following lemma.

\vspace{.1in}
\begin{lemma}
	\label{lemma: radial}
	Let $\psi: (0, \infty) \to \mathbb{R}$ be a twice differentiable function with $$\displaystyle \lim_{s \to 0} \psi (s)=0 \hspace{.4in} \text{ and} \hspace{.4in} \displaystyle \lim_{s \to 0} \psi '(s)=0. $$ Define a radially symmetric function $\omega$ on $\mathbb{R}^n$ by 
\begin{equation*}
	\omega (x) =	\left\{
	\begin{array}{l l }      
	\psi(|x|) & \hspace{.4in}  x \neq (0,...,0) , \\
		\noalign{\smallskip}
		0 & \hspace{.4in}   x = (0,...,0) .	
	\end{array}\right.
\end{equation*}
Then 
\begin{enumerate}[(a)]
	\item $\omega$ is twice differentiable everywhere in $\mathbb{R}^n$ if and only if \hspace{.02in} $\displaystyle \lim _{s \to 0}  \frac{\psi ' (s)}{s} $ exists and is finite.

	\vspace{.1in}
	
	\item $\displaystyle \lim _{|x| \to 0}\Delta \omega$ exists and is finite if and only if \hspace{.02in} $\displaystyle \lim _{s \to 0}  \psi '' (s) $ exists and is finite.
	
		\vspace{.1in}
		
	 \item  $\omega \in C^2 \left ( \mathbb{R}^n \right )$ if and only if \hspace{.02in} $\displaystyle \lim _{s \to 0}  \psi '' (s)  $ exists and  is finite.
	   
	   \vspace{.1in}
	   
	   \item $\displaystyle \lim _{|x| \to 0}\Delta \omega$ exists and is finite if and only if $\omega \in C^2 \left ( \mathbb{R}^n \right )$.
	   
	   \vspace{.1in}
	   
	   \item  $\omega$ is twice differentiable with bounded but discontinuous Hessian (and in particular bounded Laplacian) if and only if \hspace{.02in} $\displaystyle \lim _{s \to 0} \frac{\psi ' (s)}{s} $  is finite and $\psi ''(s)$ is bounded but $\displaystyle \lim _{s \to 0}  \psi '' (s)  $ does not exist.
	   
	   	\vspace{.1in}
	   	
	   \item There is no twice differentiable radial function with bounded Laplacian and unbounded Hessian.
	 
\end{enumerate}

\end{lemma}

	\vspace{.05in}

\noindent
\textbf{Remark}: By L'Hopital's Rule, if  $\displaystyle \lim _{s \to 0}  \psi '' (s)  $  is finite, so is  $\displaystyle \lim _{s \to 0}  \frac{\psi ' (s)}{s}  $.  But it is not true vice versa. For example, if $\displaystyle \psi (s)= s^4 \sin \left ( \frac{1}{s}\right )$, then $\displaystyle \psi'(s)= 4s^3\sin \left ( \frac{1}{s}\right ) -s^2 \cos \left ( \frac{1}{s}\right )  $ and $\displaystyle \psi''(s)= 12s^2\sin \left ( \frac{1}{s}\right ) -6s\cos \left ( \frac{1}{s}\right ) -\sin \left ( \frac{1}{s}\right ) $, so  $\displaystyle \lim _{s \to 0}  \frac{\psi ' (s)}{s} =0$, $\psi''(s)$ is bounded, but  $\displaystyle \lim _{s \to 0}  \psi '' (s)$ does not exist.

\vspace{.1in}
\noindent
\textit{\textbf{Proof of Lemma \ref{lemma: radial}}}:

\vspace{.05in}
By definition we only need to study the regularity of $\omega$ at the origin.  Since $\displaystyle \lim_{s \to 0} \psi(s)=0, $ $\omega$ is continuous at the origin. 
When $x \neq (0,..., 0)$, for any $i=1,...,n$,
$$ 	\frac{\partial \omega}{\partial x_i}  =  \psi ' (|x|)\frac{x_i}{|x|} \to 0  \hspace{.2in} \text{as} \,\, |x| \to 0, \hspace{.2in} \text{since} \hspace{.2in} \lim_{s \to 0} \psi '(s)=0 . $$ 

\noindent
We compute the derivative of $\omega$ at the origin by definition: 
$$\frac{\partial \omega}{\partial x_i} (0)  =  \lim _{x_i \to 0}  \frac{\omega (0, ..., x_i, ...0) - \omega(0,..., 0)}{x_i} =  \lim _{x_i \to 0}  \frac{\psi (|x_i|)}{x_i} =0, $$
where we have used $\displaystyle \lim_{s \to 0} \psi '(s)=0$ again.  Thus $\omega$ is $C^1$.  Next, we look at its second derivatives at the origin.  For any $i,j=1,...,n$ and $i \neq j$,
\begin{equation}
	\label{eqn:partial_i_j_origin}
		\frac{\partial ^2 \omega}{\partial x_i \partial x_j} (0)  =  \lim _{x_j \to 0}  \frac{\frac{\partial \omega}{\partial x_i} (0, ..., x_j, ...0) - \frac{\partial \omega}{\partial x_i}(0,..., 0)}{x_j}  =   \lim _{x_i \to 0}  \frac{0 - 0}{x_j}  =  0 .
\end{equation}
But 
\begin{equation}
	\label{eqn:partial_i_i_origin}
	\frac{\partial ^2 \omega}{\partial x_i ^2} (0)  =  \lim _{x_i \to 0}  \frac{\frac{\partial \omega}{\partial x_i} (0, ..., x_i, ...0) - \frac{\partial \omega}{\partial x_i}(0,..., 0)}{x_i} 
	=   \lim _{x_i \to 0}  \frac{\psi ' (|x_i|)}{|x_i|},
\end{equation}
 so $\displaystyle \frac{\partial ^2 \omega}{\partial x_i ^2} (0) $  is a finite number if and only if $\displaystyle \lim _{s \to 0}  \frac{\psi ' (s)}{s}  $  is a finite number.  

\vspace{.1in}
First we prove (a).
Assume $\displaystyle \lim _{s \to 0}  \frac{\psi ' (s)}{s} =L $  is a finite number, then $\frac{\partial ^2 \omega}{\partial x_i ^2 }(0) =L $ by (\ref{eqn:partial_i_i_origin}),   and

\begin{eqnarray*}
	& & \displaystyle \lim _{|x| \to 0}\frac { 	\omega (x) - \omega (0) - \displaystyle \sum _{i =1}^n \frac{\partial \omega}{\partial x_i} (0) x_i  -  
		\frac{1}{2}	\sum _{i,j=1}^n	\frac{\partial ^2 \omega}{\partial x_i \partial x_j}(0) x_i x_j  } {|x|^2} \\
	& = & \displaystyle \lim _{|x| \to 0} \frac { 	\omega (x) -   \displaystyle 
		\frac{1}{2}	\sum _{i=1}^n	\frac{\partial ^2 \omega}{\partial x_i ^2 }(0) x_i ^2  } {|x|^2} \\
	& = & \displaystyle \lim _{|x| \to 0} \frac{ \psi (|x|) -  \frac{L}{2} |x|^2 }{|x|^2} \\
	& = & 0,
\end{eqnarray*}
because by L'Hopital's Rule, $\displaystyle \lim _{s \to 0} \frac{\psi (s)}{s^2} = \lim _{x_i \to 0}  \frac{\psi ' (s)}{2s}= \frac{L}{2}$, so 
$\displaystyle \lim _{|x| \to 0}  \frac{\psi (|x|) -  \frac{L}{2} |x|^2 }{|x|^2} =0.$ Therefore, $\omega$ is twice differentiable at the origin.  

\vspace{.05in}
\noindent
On the other hand, if $\omega$ is twice differentiable at the origin, then all its first and second order partial derivatives at the origin are well defined.  In particular, $	\frac{\partial ^2 \omega}{\partial x_i ^2 }(0), i=1,...,n,$ are all defined, then $\displaystyle \lim _{s \to 0}  \frac{\psi ' (s)}{s}  $ must be finite by (\ref{eqn:partial_i_i_origin}). This proves (a).

\vspace{.1in}
Now we prove (b). First, suppose $\displaystyle \lim _{|x| \to 0}\Delta \omega$ is finite. When $x \neq (0,..., 0)$,
\begin{equation}
	\label{eqn:delta_omega}
	\Delta \omega (x)  =  \psi ''(|x|) + (n-1) \frac{\psi ' (|x|)}{|x|}.
\end{equation}
\noindent

\noindent  
Thus \, $\displaystyle \lim _{s \to 0}  \frac{s \psi ''(s) + (n-1)\psi '(s)}{s} $ \, is finite, and consequently
$$ \lim_{s \to 0} \dfrac{\frac{d}{ds} \left ( s^{n-1} \psi '(s) \right ) }{\frac{d}{ds}\left(s^{n} \right ) } = \lim_{s \to 0} \frac{s^{n-1} \psi ''(s)+ (n-1)s^{n-2}\psi'(s)}{ns^{n-1}} = \lim_{s \to 0} \dfrac{ s\psi ''(s) + (n-1) \psi ' (s)} {ns} $$

\noindent
 is finite.  Then by L'Hopital's Rule,
 \begin{eqnarray*}
 	\lim _{s \to 0} \frac{\psi ' (s)}{s} = \lim _{s \to 0} \frac{s^{n-1}\psi ' (s)}{s^n} = \lim_{s \to 0} \dfrac{\frac{d}{ds} \left ( s^{n-1} \psi '(s) \right ) }{\frac{d}{ds}\left(s^{n} \right ) }
 \end{eqnarray*}
 is finite, which implies $\displaystyle \lim_{s \to 0} \psi ''(s)$ is finite by (\ref{eqn:delta_omega}).

\noindent
Conversely, suppose $\displaystyle \lim _{s \to 0}  \psi '' (s) $ is a finite number, then $\displaystyle \lim _{s \to 0}  \frac{\psi ' (s)}{s}$ is equal to the same number.  By (\ref{eqn:delta_omega}), this implies 
$$ \lim _{|x| \to 0} 	\Delta \omega (x) = n  \lim _{|x| \to 0} \psi ''(|x|)   $$ is finite. ( In fact, by (\ref{eqn:partial_i_i_origin}) we know $\Delta \omega (0)$ is well defined and $\displaystyle  \lim _{|x| \to 0} 	\Delta \omega (x) = \Delta \omega (0).  $ )
Therefore (b) is also true.

\vspace{.1in}
Next we prove (c). 
Suppose $\displaystyle \lim _{s \to 0}  \psi '' (s) =L $  is finite, then we also know $\displaystyle \lim _{s \to 0} \frac{\psi ' (s)}{s} =L,$ and by (\ref{eqn:partial_i_i_origin}) that implies $ \frac{\partial ^2 \omega}{\partial x_i ^2} (0) =L.$   When $i \neq j$ and $x \neq (0,...,0)$, 

\begin{equation}
	\label{eqn: partial_i_j_omega}
	\frac{\partial^2 \omega}{\partial x_i \partial x_j} (x)  =   \psi '' (|x|)\frac{x_i x_j}{|x|^2} - \psi '(|x|) \frac{x_i x_j}{|x|^3}
\end{equation}
and
\begin{equation}
	\label{eqn: partial_i_i_omega}
		\frac{\partial^2 \omega}{\partial x_i^2} (x)   =	  \frac{\psi ' (|x|)}{|x|}   + \frac{x_i^2}{|x|^2}  \left ( \psi ''(|x|) - \frac{\psi '(|x|)}{|x|}   \right )	.	
\end{equation}

\noindent	
Consequently,
$$ \lim _{|x| \to 0} 	\frac{\partial ^2 \omega}{\partial x_i \partial x_j} (x)  =  \lim _{|x| \to 0} \frac{x_i x_j}{|x|^2}  \left ( \psi ''(|x|) - \frac{\psi '(|x|)}{|x|}   \right ) = 0 = \frac{\partial ^2 \omega}{\partial x_i \partial x_j} (0) . $$ Similarly,
\begin{eqnarray*}
\lim _{|x| \to 0} 	\frac{\partial^2 \omega}{\partial x_i^2} (x)  & = & 	\lim _{|x| \to 0}   \frac{\psi ' (|x|)}{|x|}   + \lim _{|x| \to 0} \frac{x_i^2}{|x|^2}  \left ( \psi ''(|x|) - \frac{\psi '(|x|)}{|x|}   \right ) \\
 & = & 	\lim _{|x| \to 0}   \frac{\psi ' (|x|)}{|x|} \\
& = & L \\
& = & \frac{\partial ^2 \omega}{\partial x_i ^2} (0) .
\end{eqnarray*}

\noindent
Therefore, $D^2 \omega$ is continuous at the origin.  Since by definition $D^2 \omega$ is also continuous away from the origin, $\omega$ is $C^2$.

\vspace{.05in}
\noindent
Conversely, suppose $\omega$ is $C^2$, then in particular $\Delta \omega$ is continuous, and by (b) $\displaystyle \lim _{s \to 0}  \psi '' (s)  $  is finite.  Thus (c) is proved, and by (b) and (c) we know (d) is true.

\vspace{.05in}
Next we look at (e).  Suppose \hspace{.02in} $\displaystyle \lim _{s \to 0} \frac{\psi ' (s)}{s} $  is finite and $\psi ''(s)$ is bounded but $\displaystyle \lim _{s \to 0}  \psi '' (s)  $ does not exist, then $\omega$ cannot be $C^2$ by (c).  By (a) $\displaystyle \lim _{s \to 0} \frac{\psi ' (s)}{s} $ being finite implies that $\omega$ is twice differentiable.  Furthermore, from (\ref{eqn:partial_i_j_origin}), (\ref{eqn:partial_i_i_origin}),	(\ref{eqn:delta_omega}), (\ref{eqn: partial_i_j_omega}), and (\ref{eqn: partial_i_i_omega}) we see that $D^2 \omega$ as well as $\Delta \omega$ are bounded everywhere.

\noindent
Conversely, suppose $\omega$ is twice differentiable but not $C^2$, and suppose $D^2 \omega$ is bounded.  First of all, by (a) the twice differentiability of $\omega$ implies \hspace{.02in} $\displaystyle \lim _{s \to 0} \frac{\psi ' (s)}{s} $  is finite. Then since $\Delta \omega$ is bounded, from (\ref{eqn:delta_omega}) we know $\psi ''(s)$ must be bounded.  However, $\displaystyle \lim _{s \to 0}  \psi '' (s)  $ cannot exist because that would imply $\omega$ being $C^2$ by (c).  This proves (e).

\vspace{.1in}
Lastly, we prove (f).  Suppose $\omega$ is twice differentiable and $\Delta \omega$ is bounded, we will show that $D^2 u$ must be bounded as well. By (a) the twice differentiability of $\omega$ implies $\displaystyle \lim _{s \to 0} \frac{\psi ' (s)}{s} $ must be finite. As a result, by (\ref{eqn:partial_i_i_origin}) we know $ \frac{\partial ^2 \omega}{\partial x_i^2}(0)$ is well defined, and by (\ref{eqn:delta_omega}) we know $\psi ''(|x|)$ must be bounded.  Then by (\ref{eqn: partial_i_j_omega}) and (\ref{eqn: partial_i_i_omega}) we know $D^2u$ is bounded for all $x \neq (0,...,0)$. Since $\frac{\partial ^2 \omega}{\partial x_i \partial x_j} (0)=0 $ when $i \neq j$, it follows that $D^2 u$ is bounded everywhere. 

\vspace{.05in}
This completes the proof of Lemma \ref{lemma: radial}.

\stop

\vspace{.4in}

\bibliographystyle{plain}
\bibliography{thesis}

\end{document}